\magnification=1200
\font\titlefont=cmcsc10 at 12pt
\hyphenation{moduli}
\input epsf

\catcode`\@=11
\font\tenmsa=msam10
\font\sevenmsa=msam7
\font\fivemsa=msam5
\font\tenmsb=msbm10
\font\sevenmsb=msbm7
\font\fivemsb=msbm5
\newfam\msafam
\newfam\msbfam
\textfont\msafam=\tenmsa  \scriptfont\msafam=\sevenmsa
  \scriptscriptfont\msafam=\fivemsa
\textfont\msbfam=\tenmsb  \scriptfont\msbfam=\sevenmsb
  \scriptscriptfont\msbfam=\fivemsb
\def\hexnumber@#1{\ifcase#1 0\or1\or2\or3\or4\or5\or6\or7\or8\or9\or
      A\or B\or C\or D\or E\or F\fi }

\font\teneuf=eufm10
\font\seveneuf=eufm7
\font\fiveeuf=eufm5
\newfam\euffam
\textfont\euffam=\teneuf
\scriptfont\euffam=\seveneuf
\scriptscriptfont\euffam=\fiveeuf
\def\frak{\ifmmode\let\next\frak@\else
 \def\next{\errmessage{Use \string\frak\space only in math mode}}\fi\next}
\def\goth{\ifmmode\let\next\frak@\else
 \def\next{\errmessage{Use \string\goth\space only in math mode}}\fi\next}
\def\frak@#1{{\frak@@{#1}}}
\def\frak@@#1{\fam\euffam#1}

\edef\msa@{\hexnumber@\msafam}
\edef\msb@{\hexnumber@\msbfam}
\mathchardef\square="0\msa@03
\mathchardef\subsetneq="3\msb@28
\mathchardef\ltimes="2\msb@6E
\mathchardef\rtimes="2\msb@6F
\def\Bbb{\ifmmode\let\next\Bbb@\else
\def\next{\errmessage{Use \string\Bbb\space only in math mode}}\fi\next}
\def\Bbb@#1{{\Bbb@@{#1}}}
\def\Bbb@@#1{\fam\msbfam#1}
\catcode`\@=12

\def\Z{{\Bbb Z}}

\def\C{{\Bbb C}}
\def\E{{\Bbb E}}
\def\F{{\Bbb F}}

\def\Q{{\Bbb Q}}
\def\PP{{\Bbb P}}
\def\F{{\Bbb F}}

\def\tto{\longrightarrow}

\def\mbar#1{{\overline{\cal M}}_{#1}}
\def\mbarg{{\overline{\cal M}}_{g}}
\def\abar#1{{\overline{\cal A}}_{#1}}
\def\Mg{{\cal M}_g}
\def\Hg{{\cal H}_g}
\def\wtmg{\widetilde{\cal M}_{g}}

\def\wthg{\widetilde{\cal H}_{g}}
\def\wtM#1{\widetilde{\cal M}_{#1}}

\def\wtH#1{\widetilde{\cal H}_{#1}}

\def\ag{\a{g}}
\def\la#1{\lambda_{#1}}

\def\m#1{{\cal M}_{#1}}

\def\a#1{{\cal A}_{#1}}

\def\ch#1{{\rm ch}_{#1}}
\def\tch{{\rm ch}}
\def\Aut{{\rm Aut}}

\def\tto{\longrightarrow}

\def\sqr#1#2{{\vcenter{\vbox{\hrule height.#2pt \hbox{\vrule width.#2pt
height #1pt \kern #1pt \vrule width.#2pt}\hrule height.#2pt}}}}
\def\qed{\hfill$\  \sqr74$}
\def\today{\ifcase\month\or
 Jan\or Febr\or  Mar\or  Apr\or May\or Jun\or  Jul\or
 Aug\or  Sep\or  Oct\or Nov\or  Dec\or\fi
 \space\number\day, \number\year}

\noindent
\font\eighteenbf=cmbx10 scaled\magstep2
\vskip 2.0pc
\centerline{\eighteenbf Complete Subvarieties of Moduli 
Spaces}
\smallskip
\medskip
\centerline{\eighteenbf and the Prym Map }
\noindent
\vskip 2pc
\font\titlefont=cmcsc10 at 11pt
\centerline{\titlefont Carel Faber and Gerard van der Geer}
\vskip 2.0pc

\bigskip
\centerline{\bf Introduction}
\medskip
In this paper we present a formula for the number of isomorphism classes
of $p$-rank zero \'etale double covers of genus $2$ curves 
over an algebraically closed field of characteristic $p>2$. 
The formula is a byproduct of our search for
complete subvarieties of moduli spaces of curves. 
Many moduli spaces are not complete because the objects that they 
parametrize can degenerate. Examples are the moduli spaces of 
polarized abelian varieties and the moduli spaces of curves. For a 
given quasi-projective variety knowledge of the maximum dimension 
of a complete subvariety gives information on its geometry and 
cohomology and it is a kind of 
measure how far the variety in question is from being affine. 
Moreover, in the case of a moduli space
it tells us when complete families must degenerate.

For the moduli space $\Mg\otimes k$ of curves of genus $g\ge2$
over a field $k$ one has 
the upper bound $g-2$ for the dimension of a complete subvariety 
due to Diaz in characteristic $0$ and Looijenga in general,
but it is not known how good this bound is. For the 
moduli space ${\cal A}_g\otimes k$ of principally polarized abelian 
varieties  of dimension $g$ over 
$k$ one has an upper bound for the dimension 
of a complete subvariety ($\dim \leq g(g-1)/2$) and it is known that 
this bound  is sharp if the characteristic of $k$ is positive. 
Recently, Keel and Sadun [KS,~Cor.1.2] proved that in char.~$0$ the bound
is {\sl not\/} sharp for $g\ge3$. In characteristic 
$p>0$ the bound is attained by the complete subvariety of abelian varieties 
of $p$-rank $0$. This suggests to look at the moduli space of curves 
in positive characteristic  and to impose conditions on the $p$-rank 
of their Jacobians in order to find complete subvarieties. 
In characteristic $p>2$ we get a subvariety of ${\cal M}_g\otimes k$ by 
considering curves $C$ of genus $g$ together with an \'etale double
cover $C'$ whose $p$-rank is $0$; its expected dimension is $g-2$.
Unfortunately, this doesn't yield a complete subvariety in general,
because there are components 
that are not complete. We hope that for almost all primes $p$ there 
are complete components and then these are necessarily 
of dimension $g-2$ and show that Diaz's bound is sharp. 
This necessitates a careful analysis of the components which 
we haven't carried out yet.
We did the analysis for genus $2$, where we have a formula for the 
number of pairs $(C', C)$ for which the $p$-rank of $C'$ is zero.
The number is zero for $p=3$, but positive for $p\geq 5$.
One of our main tools is the use of tautological classes.

The contents of this paper are as follows. 
We prove that in characteristic $p>0$ the locus of stable curves of 
$p$-rank $\leq f$ is pure of codimension $g-f$ in $\mbarg\otimes k$. 
Then we consider the Prym map and analyze it using tautological classes.
We study the locus of curves with an \'etale double cover of $p$-rank $0$
in some detail. In particular, in genus $2$ we obtain a
formula for the number of such curves. We end with several
examples illustrating our formula.

We shall use the following notation throughout this paper.
For nonnegative integers $g$ and $n$ such that $2g-2+n>0$,
we denote by ${\cal M}_{g,n}$ the moduli space of $n$-pointed curves of 
genus $g$, by $\mbar{g,n}$ the Deligne-Mumford compactification of 
${\cal M}_{g,n}$, and by $\wtM{g,n}$ the moduli space of stable
$n$-pointed curves of genus $g$ whose dual graph is a tree; 
often these curves are referred to as curves of compact type.
All these moduli spaces are smooth algebraic stacks defined over $\Z$. 
We recall that on an $n$-pointed curve the $n$ points are distinct, 
ordered, nonsingular points.
We can view $\wtM{g,n}$ as an open substack of $\mbar{g,n}$ and
${\cal M}_{g,n}$ as one of $\wtM{g,n}$. 
When $n=0$ (thus $g\ge2$), we write $\Mg$, $\mbarg$, and $\wtmg$.
We have $\mbarg = \Mg \cup \cup_{i=0}^{[g/2]} \Delta_i\,$, where $\Delta_0$
denotes the divisor of irreducible nodal curves and their
degenerations and $\Delta_i$
for $i>0$ denotes the divisor of reducible curves with components
of genus $i$ and $g-i$, and their
degenerations. Moreover, we have $\mbarg -\wtmg =\Delta_0\,$.
By $\ag$ we denote the moduli space of principally polarized abelian
varieties of dimension $g$ and by $\ag^*$ the Satake compactification
constructed by Faltings. 

In the following we shall work over an algebraically closed field $k$ 
and we shall often write $\Mg$ ($\mbarg$, $\wtmg$, etc.) instead of 
$\Mg \otimes k$ ($\mbarg \otimes k$, $\wtmg \otimes k$, etc.).
Since the stacks $\mbar{g,n}$ and $\ag$ are smooth, the codimension
of an intersection is at most the sum of the codimensions. We use this
fact without comment in \S2 and \S4.

If $k$ is of characteristic $\neq 2$ then  we denote by $R\Mg$
(resp.~$R\wtmg$)
the moduli space of pairs $(C', C)$ of connected \'etale double covers
$C' \to C$ of a genus $g$ curve $C$, with $C$ nonsingular
(resp.~of compact type). 
The curves $C'$ have genus $2g-1$. Such covers correspond exactly
to points of order $2$ in the Jacobian of $C$.

\bigskip
\centerline{\bf \S1. Complete subvarieties and tautological classes}
\medskip
It is a well-known result of Diaz [D1,~Thm.4] that a complete subvariety of the
moduli space $\Mg\otimes k$  of curves of genus $g \geq 2$ has dimension
$\leq g-2$ (the extension to positive characteristic is due to
Looijenga [L,~p.412]).
But for $g\ge4$ it is not known whether this bound is sharp.
A similar question occurs for the moduli space $\ag\otimes k$ of
principally polarized abelian varieties of dimension $g$. 
It is known by [G,~Cor.2.7]
that a complete subvariety of $\ag\otimes k$ has codimension at least
$g$ (we recall the proof below).
However, here one knows that in characteristic $p>0$ the moduli
space $\ag\otimes \F_p$ contains a complete codimension $g$ subvariety,
namely the locus of abelian varieties of $p$-rank zero, cf.~[Ko,~Thm.7]. 
As to characteristic $0$, Keel and Sadun [KS,~Cor.1.2] recently showed
that $\ag\otimes \C$ does {\sl not\/} contain a codimension $g$ 
complete subvariety for $g\ge3$, as had been conjectured by Oort [O3,~2.3G].
Thus the maximum dimension of a complete subvariety of $\ag\otimes k$
depends on the characteristic of $k\,$!
(It is not known what the maximum dimension of a
complete subvariety of $\ag\otimes \C$ is for $g>3$.)

Recall the definition of
the {\sl tautological ring}
associated to $\ag\,$. The moduli space
$\ag$ carries a tautological bundle, the Hodge bundle $\E$.
It is defined by specifying for every principally polarized
abelian scheme 
$X$ over $S$
with zero-section
$s$  a rank $g$ vector bundle on $S$: $s^*\Omega^1_{X/S}\,$. The Chern
classes $\lambda_i=c_i(\E)$ of this bundle are called {\sl tautological
classes} and they generate a $\Q$-subalgebra $T^*(\ag)$ of the Chow ring
$A^*(\ag)\otimes \Q$, called the tautological subring.

It is shown in [G,~Thm.2.1] that the relations $\ch{2k}(\E)=0$ hold $\forall
k\ge1$ in $T^*(\a{g})$, hence $T^*(\a{g})$ is a quotient of $T_g\,$,
the cohomology ring of the compact dual of the Siegel upper half space:
$$
T_g=\Q[\la1,\dots,\la{g}]/(\ch2(\E),\dots,\ch{2g}(\E))\,.
$$
Note that the vanishing of the first $g$ even $\ch{2k}$ implies the vanishing
of {\sl all\/} even $\ch{2k}$ (and in fact also the vanishing of the odd
$\ch{2k+1}$ for all $k\ge g$).

Hence $T_g$ is a complete intersection ring with socle in degree
${{g+1}\choose2}$. In [G,~Prop.2.2] it is proved that $\la{g}=0$ in
$T^*(\a{g})$, so $T^*(\a{g})$ is a quotient of $T_{g-1}\,$.
It then follows that there are no further relations: 
combine the fact
that the $p$-rank zero locus $V_0$ in $\a{g}\otimes{\overline{\F}}_p$ is
a complete subvariety of pure codimension $g$
(dimension ${g\choose2}$) 
with the ampleness of $\la1$ on
$\a{g}\otimes{\overline{\F}}_p$ due to Moret-Bailly [M-B,~p.181] 
(so $\la1^{{g\choose2}}\neq0$ in
every characteristic, and a generator of the socle of $T_{g-1}$
doesn't vanish on $\a{g}$). To summarize the above:
$$
T^*(\a{g})\cong T_{g-1}\,.
$$
This implies that  $\la1^{{g \choose 2} +1}=0$ and hence that a complete
subvariety of $\ag$ has codimension at least $g$.

For a more complete discussion see [FL,~\S11] and [GO,~\S3], but note that 
Thm.~3.4 in [GO] is at present only a conjecture, unfortunately.

\bigskip
\centerline{\bf \S2. Subvarieties of the moduli space of curves defined by
the $p$-rank}
\medskip
We are interested in the loci of curves whose $p$-rank is $0$.
For a smooth curve $C$ of genus $g$ over an algebraically closed field $k$
of characteristic $p>0$ 
we define the $p$-rank $f$ of $C$ as the $p$-rank of its Jacobian $J$:
if $J[p]$ denotes the kernel of multiplication by $p$ on $J$ then
$$
\# J[p](k)= p^f.
$$
We have $0\leq f \leq g$. Alternative definitions of this $p$-rank
are: the semisimple rank of the $\sigma$-linear operator $F$ induced
by the relative Frobenius on $H^1(J,O_J)$ (or on $H^1(C,O_C)$); or
dually, the semisimple rank of the $\sigma^{-1}$-linear Cartier operator $V$
on $H^0(J,\Omega^1_J)$ (or on $H^0(C,\Omega_C^1)$).

To extend this definition to stable curves we consider a family $C
\to B$ of stable curves of genus $g$ over an irreducible base scheme
$B$. We let $\omega_{C/B}$ be the relative dualizing sheaf and write 
$\E=\pi_*(\omega_{C/B})$ for the Hodge bundle, a locally free
sheaf of rank $g$. On it we have the relative Cartier operator $V:
\E \to \E$. For a point $b \in B$ the $p$-rank of the fibre $C_b$ is
defined as the semisimple rank of $V$ on the fibre $\E_b$.
Equivalently, we can look at the action of the relative Frobenius on
$R^1\pi_*O_C$. A result of Grothendieck (see [Ka,~Th.2.3.1]) says that the
$p$-rank is lower semi-continuous on $B$.

The following lemma is a variant of a
lemma of Oort, see [O1,~1.6]. The proof is similar to that of Oort
and is omitted.

\proclaim (2.1) Lemma. Let $B$ be an irreducible scheme over $k$ and $C\to
B$ a stable curve over $B$. Let $f$ be the $p$-rank of the generic
fibre $C_{\eta}$. Let $W\subset B$ denote the closed set over which the
$p$-rank of the fibre is at most $f-1$. Then either $W$ is empty or
$W$ is pure of codimension $1$ in $B$.\qed \par

\noindent
{\bf (2.2) Definition.} We define the locus $V_f(\mbarg)$  as
the locus of curves with $p$-rank $\leq f$ in $\mbarg$~;
similarly, we define $V_f(R\wtmg)$ as the locus of pairs $(C', C)$
for which the
$p$-rank of $C'$ is $\leq f$.  We denote  by $RV_0(\wtmg)$ the
locus of pairs $(C', C)$ for which the $p$-rank of $C$ is $0$.
Similarly, we denote by  $V_f(\ag)$ the locus of principally polarized
abelian varieties of dimension $g$ with $p$-rank $\leq f$.

\medskip
We have the following result.

\proclaim (2.3) Theorem. The locus $V_f({\mbarg})$ is pure of
codimension $g-f$ in $\mbarg$~.
\par
\noindent
{\sl Proof.} 
We apply Lemma 2.1 just as in [Ko,~p.164]: 
Fix $r<g$. Let $C_r$ be any component of $V_r(\mbarg)$. 
Let $C_{r+1}$ be a component of $V_{r+1}(\mbarg)$ containing $C_r$. 
Note that {\sl a priori\/} it is possible that $C_r = C_{r+1}$. We obtain
a sequence $C_r\subseteq C_{r+1} \subseteq \dots \subseteq C_g = \mbarg$. 
The lemma tells us: for $r'=r,r+1,\dots,g-1$, if $C_{r'}\neq C_{r'+1}$,
then $C_{r'}$ has codimension one in $C_{r'+1}$. We conclude: the codimension
of $C_r$ is at most $g-r$.

Now let $r=0$. As is well-known, $C_0$ is a complete subvariety of
$\wtmg$ (the $p$-rank of a generalized abelian variety is at least its
torus rank). 
On the other hand, we have Diaz's upper bound $2g-3$ for the
dimension of such a complete subvariety; see Lemma 2.4 below.
It follows that $C_0$ has codimension $g$: the case $f=0$ of the theorem.

If every $C_r$ would contain a $C_0$, we would now be done.
Since we don't know how to establish this, we proceed differently.
We prove that every $C_r$ has codimension $g-r$ in $\mbarg$ by induction
on $r$. The case $r=0$ has been established, so assume $r>0$.
Since ${\rm codim}(C_r)\leq g-r\leq g-1$, the component $C_r$
cannot be complete in $\wtmg$,
hence intersects $\Delta_0$. The intersection consists of points corresponding
to curves of genus $g-1$ and $p$-rank $\leq r-1$ with two points identified.
The dimension of the intersection is at most $2(g-1)-3+(r-1)+2=2g-4+r$:
this is obvious for $g=2$ and follows for $g\ge3$ by induction.
Hence the dimension of $C_r$ is $2g-3+r$. \qed

\medskip \noindent {\bf Remark.} 
It is not clear whether $V_f(\ag)$ intersects the Torelli image
of $\wtmg$ transversally. 
Also, is for example $V_0(\mbarg)$ irreducible for $g\ge4$~?

\medskip
The theorem we just established is the analogue for the moduli space of
curves of a result of Koblitz [Ko,~Thm.~7] for the moduli of principally 
polarized abelian varieties and of the main result of
Norman and Oort [NO,~Thm.4.1] for abelian varieties with 
arbitrary polarization. 


The following result was already observed by Diaz [D2,~p.80].

\proclaim (2.4) Lemma. The dimension of a complete subvariety of
$\wtM{g,n}$ is at most $2g-3+n$. \par

\noindent{\sl Proof.} This is trivial for $g\le1$. For $g\ge2$, it
suffices to prove this for $n=0$. By the Diaz-Looijenga bound $g-2$
for the dimension of a complete subvariety of $\Mg$~, we may assume
that the complete subvariety $Y$ meets some $\Delta_i$, with $0<i\le g/2$.
Then by induction $\dim Y\le 1+(2i-2)+(2(g-i)-2)=2g-3$. \qed

\proclaim (2.5) Lemma. Let $X$ be a complete subvariety of $\wtmg$
of dimension $2g-3$. Then $X$ contains points corresponding to chains
of $g$ elliptic curves. In particular, $X$ intersects every component
of the boundary $\wtmg-\Mg$.

\noindent
{\sl Proof.} 
Since $X$ is of maximal dimension, the inequalities in the proof 
of Lemma 2.4 are equalities. Some geometric consequences of this will
be exploited here. We prove the statement by induction on $g$; the case
$g=2$ is clear, so assume $g\ge3$.
We know that $X$ meets $\Delta_i$ for some $i$ with $0<i\le g/2$.
The inverse image
of $X\cap\Delta_i$ in $\wtM{i,1}\times\wtM{g-i,1}$ is pure of dimension $2g-4$.
Let $Z$ be a component. Let $A_{i,1}$ be the projection on the first factor
and let $B_{g-i,1}$ be the projection on the second factor. Then
$\dim A_{i,1}=2i-2$ and $\dim B_{g-i,1}=2(g-i)-2$. It follows that
$Z=A_{i,1}\times B_{g-i,1}$.
When $i=1$, the factor $A_{1,1}$ corresponds to a single elliptic curve.
The image $B_{g-1}$ in $\wtM{g-1}$ and, for $i\ge2$, 
the images $A_i$ and $B_{g-i}$ in $\wtM{i}$ resp.~$\wtM{g-i}$
contain points corresponding to smooth curves; this follows again from
maximality. Hence $A_{i,1}$ is the inverse image of $A_i$ for $i\ge2$
and $B_{g-i,1}$ is the inverse image of $B_{g-i}$ for all $i$.
Now use induction: pick points in $A_i$ (for $i\ge2$) and $B_{g-i}$
corresponding to chains of elliptic curves, and (for any $i$) pick points
in $A_{i,1}$ and $B_{g-i,1}$ corresponding to marked points on 
extremal components of the chains. Gluing gives a chain of $g$ elliptic
curves. \qed

\medskip \noindent {\bf Remark.} 
The proof above is very similar to the way Keel and Sadun deduce from 
[KS,~Cor.1.2] that a complete $X$ of dimension $2g-3$ in $\wtmg$
doesn't exist in characteristic $0$ for $g\ge3$.

\medskip
Let $\Hg$ be the hyperelliptic locus in $\Mg$ and $\wthg$
the hyperelliptic locus in $\wtmg$, that is, the closure 
of $\Hg$ in $\wtmg$. In what follows, we assume that the characteristic
is different from $2$. Then the theory of admissible covers [HM,~\S4]
can be used to describe the stable hyperelliptic curves.

For $g\ge1$, let $(D,b_1,\dots,b_{2g+2})$ be a stable $(2g+2)$-pointed
curve of genus $0$. The dual graph of $D$ is a tree.
Define a node of $D$ to be even (resp.~odd) if the
number of marked points on either side of the node is even (resp.~odd).
There exists a unique admissible cover $A$ of degree $2$ of $D$. The inverse
image of a component of $D$ is the unique double cover of that component
that is ramified exactly over the marked points and the odd nodes on that
component. (The double cover is disconnected exactly when the component
contains no marked points and only even nodes.) The $2g+2$ Weierstrass points
on $A$ are by definition the inverse images of the $2g+2$ marked points.
The admissible cover is a nodal curve. It may contain smooth rational components
that meet the rest of the curve in only $2$ points. (Observe that these
necessarily cover a component of $D$ with $2$ marked points and $1$ node.)
A stable
hyperelliptic curve of genus $g\ge2$ is by definition the stable curve of genus
$g$ obtained from the admissible cover by contracting those components.
The Weierstrass points on a stable hyperelliptic curve are the images of the
Weierstrass points on the admissible cover. Only the Weierstrass points
on contracted components become singular points.
When $g=1$, we remember the Weierstrass points and consider the admissible
covers as stable $4$-pointed curves of genus $1$.

A stable hyperelliptic curve is of compact type if and only if all the nodes
of the stable $(2g+2)$-pointed curve $D$ of genus $0$ are odd. No components
are contracted and these curves may be identified with the admissible covers.
Therefore $\wthg$, considered as a coarse moduli space, may be identified
for $g\ge1$
with the quotient of an open set in $\mbar{0,2g+2}$ (the complement
of the `even' boundary divisors) by the natural action of the symmetric group
on $2g+2$ letters.

\proclaim (2.6) Lemma. {\rm (Char.~$\neq2$.)} A complete subvariety of
$\wthg$ has dimension at most $g-1$. If $Z$ is a complete 
subvariety of dimension $g-1$ then $Z\cap \Hg\neq \emptyset$.
\par
\noindent
{\sl Proof.}
By induction on $g$, the case $g=1$ being trivial. 
Assume $g>1$ and let $Z$ be 
a positive-dimensional complete subvariety of $\wthg$. Then
$Z\cap\Delta_i\neq\emptyset$ for some $i$ with $0<i\le g/2$, since
$\Hg$ is affine. Denote by $W\wtH{i}\subset\wtH{i,1}$ the locus where
the marked point is a Weierstrass point. There exists a finite and
surjective map $W\wtH{i}\times W\wtH{g-i}\to \Delta_i\cap\wthg$.
Hence $\dim Z\le 1+(i-1)+(g-i-1)=g-1$. If equality holds then $Z$
is not contained in the boundary. \qed

\let\xpar=\par
\medskip
In characteristic $0$ a complete subvariety of $\wtmg$ of dimension
$2g-3$ doesn't exist for $g\ge3$.
In positive characteristic we have the following.

\proclaim (2.7) Proposition. {\rm (Char.~$>2$.)}
Let $X$ be a complete subvariety of $\wtmg$ of maximal dimension $2g-3$, e.g.,
an irreducible component of $V_0(\mbarg)$.
\item{(i)} $X$ contains
points corresponding to chains of $g$ elliptic curves such that
on every non-extremal elliptic curve in the chain the two
attachment points differ by a 2-torsion point.
\item{(ii)} $X$ contains points corresponding to smooth hyperelliptic curves.

\noindent
{\sl Proof.} 
The proof of the first statement is entirely analogous to the proof
of Lemma 2.5; the only adaptation required is that one should choose
points in $A_{i,1}$ and $B_{g-i,1}$ corresponding to marked points on extremal
components of the chains that differ from the attachment point by a point
of order 2 (such points exist). Fulton ([HM,~p.88]) shows that the curves
from (i) are stable hyperelliptic curves. Thus $X\cap\wthg$ is complete of
dimension at least $g-1$ and nonempty. We conclude by applying Lemma 2.6. \qed

\bigskip
\centerline{\bf \S3. Relations between the tautological classes}
\medskip
We work over an algebraically closed field of characteristic $\neq 2$.
Let as before $R\Mg$ denote the moduli space of pairs $(C',C)$
with $C$ a smooth curve of genus $g$ and $C'$ an \'etale double
cover of $C$. We consider the Prym map
$$
P : R\Mg \tto \a{g-1}, \qquad (C',C) \mapsto P(C',C).
$$
Here $P(C',C)$ is the Prym variety of $C' \to C$, i.e.,
the identity component of the kernel of the norm map ${\rm Nm}: J(C')
\to J(C)$. The Prym variety comes with a canonical principal
polarization $\Xi$, see [M1,~p.333].

The Prym map $P$  extends  to a map
$$
P: R\wtmg \tto \a{g-1}^*,
\qquad (C',C) \mapsto P(C',C),
$$
where $\a{g-1}^*$ is the Satake compactification and
$R\wtmg$ is the moduli space of pairs $(C',C)$ with $C$
a stable curve of genus $g$ of compact type 
and $C'$ an \'etale double cover of $C$. Then $C'$ is a stable curve
of genus $2g-1$, not necessarily of compact type. The
Prym variety is then the identity component of the kernel of the norm map
on the generalized Jacobians.

We establish a relation between the Chern classes
of the Hodge bundle and the pull-back of such classes under the Prym
map. Consider the universal morphism of curves $q: C'\to C$
over the moduli space $R\wtmg$.  The universal curve $\pi:C=C_g\to R\wtmg$
carries a torsion line bundle $L$ defined by $q_*(O_{C'}) = O_C
\oplus L$ and $L^{\otimes 2} \simeq O_C$.

Let $\omega$ be the relative dualizing sheaf of
$C_g$ over $R\wtmg$.  We define the Prym-Hodge bundle $\E^{\prime}$ on 
$R\wtmg$  by
$$
\E^{\prime}:=\pi_*(\omega\otimes L).
$$
This is a vector bundle of rank $g-1$ on $R\wtmg$. We denote its
Chern classes by $\lambda_i^{\prime}=c_i(\E^{\prime})$,
for $1\le i\le g-1$.

\proclaim (3.1) Theorem. Let $\phi: R\wtmg \to \wtmg$ be the map
$(C',C)\mapsto C$. Then  $\tch(\E^{\prime})=\tch(\phi^*(\E))-1.$
\par
\noindent
{\sl Proof.} Applying the Grothendieck-Riemann-Roch theorem (GRR)
to the line bundle $\omega\otimes L$ and
the morphism $\pi$ gives
$$
\tch(\E^{\prime}) =   \pi_*(\tch(\omega\otimes L)\cdot 
{\rm Td}^{\vee}(\omega)).
$$
Since we work in Chow groups tensored with $\Q$ we can disregard
torsion classes and we may replace the term $\omega \otimes L$ on the
right by $\omega$, so that
$$
\tch(\E^{\prime})= \pi_*(\tch(\omega)\cdot {\rm Td}^{\vee}(\omega)).\eqno(1)
$$
But applying GRR to $\omega$ and $\pi$ gives
$$
\tch(\pi_! \omega)=\tch(\phi^*(\E))-1 =  \pi_*(\tch(\omega )\cdot {\rm
Td}^{\vee}(\omega)), \eqno(2)
$$
so that by comparing (1) and (2) we find $\tch(\E^{\prime})=
\tch(\phi^*(\E))-1$ as required. \qed

\proclaim (3.2) Corollary. We have $\phi^*(\lambda_i)=\lambda_i^{\prime}$
for $i=1,\ldots,g-1$.\qed

Note that it also follows that $\lambda_g$ vanishes in the Chow ring 
$A^*(\wtmg)$. This is compatible with the fact that it vanishes on $\ag$.

Denote by
$\psi: R\wtmg\to\mbar{2g-1}$ the morphism $(C',C)\mapsto C'$.
Then $\psi^*(\E_{2g-1})=\phi^*(\E)\oplus\E'$, hence $\psi^*(\la1)
=2\la1'$.

\proclaim (3.3) Corollary. 
The Torelli maps $R\wtmg\to\a{g}^*$ (resp.~$R\wtmg\to\a{2g-1}^*$)
sending $(C',C)$ to $J(C)$ (resp.~$J(C')$) are constant on any complete
connected algebraic subset of 
a fiber of the extended Prym map $P: R\wtmg\to\a{g-1}^*.$
\par
\noindent
{\sl Proof.} 
Clearly $P^*(\la1)=\la1'$. This is zero on a fiber of $P$, hence the
same holds for $\phi^*(\la1)$ and $\psi^*(\la1)$. These are pull-backs
of ample classes on $\a{g}^*$ resp.~$\a{2g-1}^*$ via the Torelli maps.
The result follows. \qed

\proclaim (3.4) Corollary. The restriction of the Prym map to
a complete subvariety $Y$ of $R\wtmg$ is quasi-finite on $Y\cap R\Mg$.
\qed

\bigskip
\centerline{\bf \S4. The locus $V_0(R\Mg)$}
\medskip
We now consider in $R\Mg$
the locus $V_0(R\Mg)$ of \'etale double covers $(C',C)$
where the $p$-rank of $C'$ is zero ($p>2$). Here is our motivation for
studying this locus. Consider $\ag$ and a toroidal compactification
$\abar{g}$ of the type constructed in [FC,~Ch.IV] (such that the
Hodge bundle extends).
The class $\la{g}$ has zero intersection with all boundary
classes (those from $\abar{g}-\ag$). If a positive multiple of $\la{g}$
can be written as an effective sum of subvarieties not contained in the
boundary, then all those subvarieties are complete subvarieties
of $\ag$ of codimension $g$.
The locus $V_0$ is a (canonical) effective cycle of this type in positive
characteristic [G,~Thm.9.2], while in characteristic $0$ such an effective cycle
doesn't exist for $g\ge3$ [KS,~Cor.1.2]. As we saw in \S2, the situation 
for $\wtmg$ and its compactification $\mbar{g}$ is exactly analogous.

Now consider $\Mg$, with compactification $\mbar{g}$. The class
$\la{g}\la{g-1}$ has zero intersection with all boundary classes [F,~p.112]
and appears to be the natural analogue of the class $\la{g}$
for $\ag$ and $\wtmg$ (cf.~[FP,~\S0]). In the search for an effective
representative of a multiple of $\la{g}\la{g-1}$, it seems natural
to combine the $p$-rank zero requirement with a condition specific
to curves. Certainly, the construction of the Prym variety
$P(C',C)$ associated to an \'etale double cover $C'$ of a curve $C$
plays an important role in the theory of curves and Jacobians. 
The Prym variety is
a principally polarized abelian variety of dimension $g(C)-1$. The
requirement that it have $p$-rank zero should then impose the desired
$g-1$ additional conditions; or equivalently, one may combine the two
requirements by asking that $C'$ have $p$-rank zero. 
Note also that the class of $V_0(\mbar{2g-1})$ is a multiple of $\la{2g-1}$
and that $\psi^*(\la{2g-1})=\phi^*(\la{g}\la{g-1})$.
This explains our interest in the locus $V_0(R\Mg)$.

\proclaim (4.1) Proposition. Every irreducible component of $V_0(R\Mg)$
has dimension at least $g-2$. A complete irreducible component of $V_0(R\Mg)$
has dimension equal to $g-2$.

\noindent {\sl Proof.} The first statement follows by applying Lemma 2.1
to the family of curves $C'$ over $R\Mg$. The second statement follows
then from Diaz-Looijenga. \qed

\medskip
Let as before
$\psi: R\wtmg\to\mbar{2g-1}$ be the morphism $(C',C)\mapsto C'$.
We define
$$V_0(R\wtmg)=\psi^{-1}(V_0(\mbar{2g-1})).$$
Clearly, this contains $V_0(R\Mg)=\psi^{-1}(V_0(\m{2g-1}))$, and every
irreducible component of $V_0(R\wtmg)$ has dimension at least $g-2$,
by Lemma 2.1 again. Observe that $V_0(R\wtmg)$ is complete: it is a 
closed sublocus of the complete locus $RV_0(\wtmg)$.

For $g=2$, all components of $V_0(R\wtmg)$ have the expected dimension
$g-2=0$. In fact, from Corollary 3.3 and the identification
$\wtM2=\a2$, the Prym map $P$ is finite on $RV_0(\wtM2)$.
In \S6 we count the number of points of $V_0(R\m2)$.
It turns out that $V_0(R\m2)$ is empty for $p=3$, but nonempty for $p\ge5$.

For $g\ge3$, the situation is more complicated. 
We first analyze the case $g=3$ in some detail.
We have the following description of the components of $V_0(R\wtM3)$.

\proclaim (4.2) Proposition. 
\item{(i)} All components of $V_0(R\m3)$ have dimension equal to $g-2=1$.
\item{(ii)} For every pair $(D',D)\in V_0(R\m2)$ and every supersingular
elliptic curve $E$, there is a component of $V_0(R\wtM3)$ isomorphic to $D$.
The curve $C_p$ corresponding to $p\in D$ is obtained by gluing $D$ and $E$
at $p$, while $C'_p$ is obtained by gluing two copies of $E$ to $D'$
at the two inverse images of $p$.
\item{(iii)} For every pair $(E',E)\in V_0(R\m{1,1})$ and every
component $X$ of $V_0(\wtM2)$, there is a component of $V_0(R\wtM3)$
isomorphic to the universal curve over $X$ (considered as a subvariety
of $\wtM{2,1}$). The curve $C$ corresponding
to $D\in X$ and $p\in D$ consists of $E$ and $D$ glued at $p$, while
$C'$ is obtained by gluing two copies of $D$ to $E'$ (the inverse images
of $0\in E$ are identified with the two points $p$).
\item{(iv)} For every pair $(E',E)\in V_0(R\m{1,1})$ and every two
supersingular elliptic curves $F$ and $G$, 
there is a component of $V_0(R\wtM3)$ isomorphic to $E$ (considered as a
subvariety of $\wtM{1,2}$). The curve
$C_p$ corresponding to $0\neq p\in E$ consists of $E$ with $F$ and $G$
glued on at $0$ and $p$, while $C'_p$ consists of $E'$ with two copies
of both $F$ and $G$ glued on at the respective inverse images.
\xpar\noindent {\sl The components of type (ii)--(iv) are the components
of $V_0(R\wtM3)$ not intersecting $R\m3$.}

\epsfbox{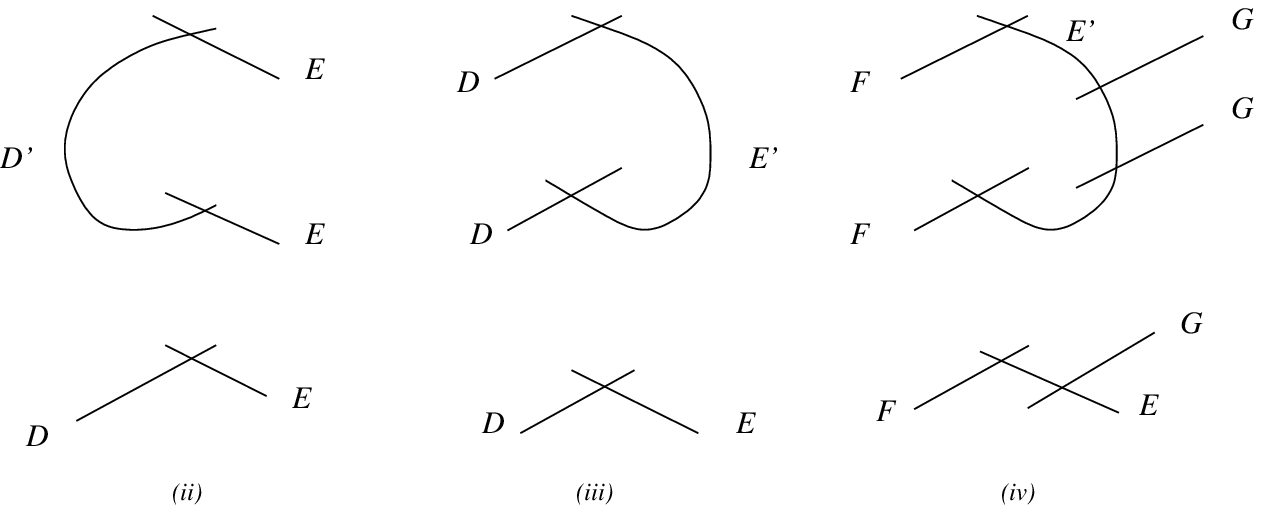}

(Figures (iii) and (iv) represent general points of the
corresponding components.)
\medskip

\noindent {\sl Proof.} 
From Corollary 3.4, the restriction of the Prym map to $RV_0(\Mg)$
is quasi-finite. For $g=3$ this means that all components of $V_0(R\m3)$
have dimension equal to $g-2=1$, the dimension of $V_0(\a2)$.
By taking the closure, we obtain
the components of $V_0(R\wtM3)$ intersecting $R\m3$.
It is straightforward to check 
that (ii), (iii), and (iv) give a complete description of the
$p$-rank zero \'etale double covers of reducible curves of genus $3$.
Since these families are positive-dimensional, they yield 
the components of $V_0(R\wtM3)$ not intersecting $R\m3$. \qed 

\medskip
The components of type (iii) have dimension $2$ and provide the first
examples of components of $V_0(R\wtmg)$ of dimension greater than $g-2$.
Taking $X$ to be a component of $V_0(\wtM{g-1})$, we find components
of $V_0(R\wtmg)$ of dimension $2g-4$.

Return to genus $3$. From the description in Proposition 4.2, one might
expect that all components of $V_0(R\m3)$ are complete. However, this
is not the case, as is shown by the following argument that we learned
from Oort.

\proclaim (4.3) Proposition. {\rm (Oort.)} If $V_0(R\m{g-1})$ is nonempty,
then it has a component of dimension $>g-3$
(necessarily noncomplete) or $V_0(R\m{g})$ has a noncomplete component.

\noindent {\sl Proof.} 
Consider $\psi:R\wtmg\to\mbar{2g-1}$ and the Torelli map
$t:\mbar{2g-1}\to\a{2g-1}^*$. The image $t\psi(R\wtmg)$ has dimension
$3g-3$; every component of the intersection with $V_0(\a{2g-1})$ has
dimension at least $g-2$. The intersection 
$V_0(\a{2g-1}) \cap t \psi (R \wtmg)$ equals $t\psi(V_0(R\wtmg))$.

Let $Y$ be a component of $V_0(R\m{g-1})$. We may assume that
$\dim Y=g-3$. Consider the following family in $V_0(R\wtmg)$ isomorphic to the
universal curve over $\phi(Y)$: for $(D',D)\in Y$ and $p\in D$, construct
$C$ by gluing a supersingular elliptic curve $E$ to $D$ at $p$
and $C'$ by gluing two copies of $E$ to $D'$ (just as in type (ii) above).
Denote this family by $Z$. Then $t\psi(Z)$ has dimension $g-3$.
Therefore it is strictly contained in a component of $t\psi(V_0(R\wtmg))$.
That component necessarily contains Jacobians of smooth curves $C'$
and yields a noncomplete component of $V_0(R\m{g})$. \qed

\medskip
We may similarly consider a family of pairs $(C',C)$, where $C$ consists
of two smooth components, $D$ of genus $i$ and $E$ of genus $j$, and $C'$
consists of $D'$ of genus $2i-1$ and two copies of $E$. 
Let $(D',D)$ vary in a component $Y$ of $V_0(R\m{i})$ and $E$ in a component
of $V_0(\m{j})$. Since the attachment points vary also, the dimension
of this family equals $(\dim Y+1)+(2j-2)\ge(i-1)+(2j-2)=j+g-3$.
({\sl Mutatis mutandis\/}, this
includes the cases $i=1$ and $j=1$.) Denoting the family by $Z$, we have
$$\dim t\psi(Z)=\left\{
\eqalign{&0,\cr
&j+g-4=2g-5,\cr
&\dim Y\ge i-2=g-3, \cr
&\dim Y+2j-3\ge j+g-5,\cr
}\right.
\qquad
\eqalign{&i=1,\quad j=1;\cr
&i=1,\quad j>1;\cr
&i>1,\quad j=1;\cr
&i>1,\quad j>1.\cr
}
$$
The proposition uses the case $j=1$ for $g\ge3$. We see that for $g\ge4$
one may also use the case $j=2$. This proves: if $V_0(R\m{g-2})$ is nonempty,
then it has a component of excess dimension 
or $V_0(R\m{g})$ has a noncomplete component.

\proclaim (4.4) Corollary. For $p\ge5$, the locus $V_0(R\Mg)$ has a
noncomplete component for $g=3$, $4$, and $5$.

\noindent {\sl Proof.} Use that $V_0(R\m2)$ is nonempty and of dimension $0$.
Then $V_0(R\m3)$
has a noncomplete component. It has dimension $1$. Now use the proposition
again for $g=4$ and the discussion above for $g=5$.\qed

\medskip
We study the case $p=3$ and $g=3$ in the next section.

\bigskip
\centerline{\bf \S5. The locus $V_0(R\m3)$ for $p=3$}
\medskip
We start by recalling some facts concerning hyperelliptic curves.
Assume ${\rm char}(k)\neq2$.
Every \'etale double cover
$C'$ of a hyperelliptic curve $C$ arises from splitting up the set $B$
of branch points of $C$ into two disjoint subsets of even cardinality:
$$B=B_1\amalg B_2. $$
Then $C'$ is a Galois cover of $\PP^1$ with Galois group
$\Z/2 \times \Z/2$ and we have a diagram
$$
\matrix{
&&& C' &&&\cr
&&\swarrow & \downarrow & \searrow&\cr
& C_1 && C && C_2 &\cr
&&\searrow & \downarrow & \swarrow && \cr
&&& \PP^1 &&&\cr
}
$$
where $C_1$ (resp.\ $C_2$) is the hyperelliptic
curve with branch points $B_1$ (resp.\ $B_2$).
Moreover, the Prym variety of $C'$ over $C$ is isomorphic to
$J(C_1)\times J(C_2)$ (cf.~[M1,~p.346]).
Note that $C'$ is hyperelliptic if and only if $C_1$ or $C_2$
is a rational curve.

Next, assume that $k$ has characteristic $p>2$. The $p$-rank of a
hyperelliptic curve $C$ can be computed easily. If $C$ is given
by $y^2=f(x)$ with $f(x)$ a polynomial 
without multiple roots of degree $2g+1$ or $2g+2$, then 
the regular differentials $x^{i-1}dx/y$ with $1\le i\le g$ form a basis.
Write 
$$
f(x)^{(p-1)/2}=a(x)= \sum_{i=0}^{\infty} a_i x^i,
$$
then the Hasse-Witt matrix with respect to this basis is 
$$
H=\left( \matrix{
a_{p-1} & a_{2p-1} & \ldots & a_{gp-1} \cr
a_{p-2} & a_{2p-2} & \ldots & a_{gp-2} \cr
\vdots & \vdots & \ddots & \vdots \cr
a_{p-g} & a_{2p-g} & \ldots & a_{gp-g} \cr}
\right).
$$
(Cf.~[SV,~p.54].) The $p$-rank of $C$ equals the semisimple rank of $H$.
In particular, $C$ has $p$-rank $0$ if and only if
$$H\cdot H^{(p)} \cdot \dots \cdot H^{(p^{g-1})}=0,$$
where $H^{(p^i)}$ is the matrix obtained from $H$ by raising every
entry to the $(p^i)$th power. This is independent of the chosen basis.

\medskip
Let now $p=3$.

\proclaim (5.1) Lemma. A hyperelliptic curve $C'$ of genus $5$ with
a fixed-point-free involution cannot have $3$-rank $0$.

\noindent {\sl Proof.} 
Let $C$ be the quotient curve of genus $3$. Then $C$ is hyperelliptic
and $C'$ is an \'etale double cover of $C$. By the discussion above,
the Prym variety of $(C',C)$ is the Jacobian of a curve $D$ of genus $2$.
Assume that $C'$ has $3$-rank $0$. Then $D$ has $3$-rank $0$. Write $D$
as $y^2=\sum_{i=0}^5\,a_ix^i$. The Hasse-Witt matrix of $D$ equals
$$H_D=\left( \matrix{a_2 & a_5 \cr a_1 & a_4 \cr}\right).$$
We may assume $a_5=1$ and $a_4=0$. It follows then that $a_1=a_2=0$.
We may assume $a_0=1$ and write $D$ as
$$y^2=x^5+ax^3+1.$$
Then $C$ can be written as
$$y^2=(x^2+bx+c)(x^5+ax^3+1)=x^7+bx^6+(a+c)x^5+abx^4+acx^3+x^2+bx+c$$
and the Hasse-Witt matrix of $C$ equals
$$H_C=\left( \matrix{ 1 & a+c & 0 \cr b & ab & 1 \cr c & ac & b }\right),$$
with determinant $c^2-cb^2=c(c-b^2)$. The discriminant of $x^2+bx+c$
equals $b^2-c$, so it follows that $c=0$. Then 
$$H_C=\left( \matrix{ 1 & a & 0 \cr b & ab & 1 \cr 0 & 0 & b }\right)$$
and the $(3,3)$-entry of $H\cdot H^{(3)}\cdot H^{(9)}$ equals $b^{13}$.
This forces $b=0$ and we obtain a contradiction. \qed

\proclaim (5.2) Proposition. For $p=3$, the locus $V_0(R\m3)$ doesn't have
a complete component.

\noindent {\sl Proof.}
Suppose instead that there is a complete one-dimensional family $X$
of smooth curves of genus $5$, varying in moduli, such that each curve
has $3$-rank $0$ and possesses a fixed-point-free involution. On $X$
the degree of $\la1$ is positive. Hence the degree of the divisor $T$
in $\m5$ of trigonal curves and their degenerations is positive,
since $[T]=8\la1$ (cf.~[HM,~p.24]).
It is easy to see that a trigonal curve of genus $5$
cannot have a fixed-point-free involution: the $g^1_3$ is unique, thus
fixed by the involution; there are $2$ fixed divisors; this forces
fixed points on the curve. The intersection points of $X$
and $T$ come therefore from degenerate trigonal curves
where the trigonal system decomposes as a $g^1_2$ plus a base point.
These curves are hyperelliptic and Lemma 5.1 gives a contradiction. 
\qed

\medskip\noindent {\bf (5.3) Example.}
For $p=5$, there do exist smooth hyperelliptic curves of genus
$5$ with $5$-rank $0$ and possessing a fixed-point-free involution. 
Consider $C$ of genus $3$ given by
$$
y^2=(x^2-x+2)(x^5+x^4+2x^2-2x)
=x^7+x^5-x^4+x^3+x^2+x.
$$
One easily checks that both $C$ and $D$, given by $y^2=f(x)=x^5+x^4+2x^2-2x$,
have $5$-rank $0$. The \'etale double cover $C'$ is given by
$$y^2=f\bigg({{x^2-2}\over{2x-1}}\bigg)(2x-1)^6=
2x^{11}-2x^{10}+x^9+x^8-x^7-2x^6+2x^5-x^4+x^3-2x^2-2x+2,$$
with fixed-point-free involution $(x,y)\mapsto({{x+1}\over{2x-1}},-y)$.



\medskip\noindent {\bf (5.4) Example.}
For $p=3$, there do exist smooth hyperelliptic curves of genus
$7$ with $3$-rank $0$ and possessing a fixed-point-free involution.
Consider $C$ of genus $4$ given by
$$
y^2=x(x^8+x^6+x^5-x^3-x^2-1).
$$
One easily checks that both $C$ and $D$, 
given by $y^2=f(x)=x^8+x^6+x^5-x^3-x^2-1$,
have $3$-rank $0$. The \'etale double cover $C'$ is given by
$y^2=f(x^2)=x^{16}+x^{12}+x^{10}-x^6-x^4-1$,
with fixed-point-free involution $(x,y)\mapsto(-x,-y)$.

\medskip\noindent {\bf (5.5) Example.}
For $p=3$, there do exist smooth hyperelliptic curves of genus
$3$ possessing an \'etale double cover with $3$-rank $0$. 
Consider $C$ of genus $3$ given by
$$y^2 = (x^3+x+1) (x^4+x-1)= x^7+x^5-x^4-x^3+x^2-1.$$
One easily checks that $C$ and the elliptic curves $v^2=x^3+x+1$
and $w^2=x^4+x-1$ have $3$-rank $0$. The \'etale double cover $C'$ is 
the normalization of
$$v^8-w^6-v^4w^2-v^6-w^4-v^4+v^2-1=0,$$
with fixed-point-free involution induced by $(v,w)\mapsto(-v,-w)$.

\bigskip
\centerline{\bf \S6. The number of $p$-rank zero \'etale double
covers of genus $2$ curves}
\medskip
In this section we count the number of points of $V_0(R\m2)$, that is,
we count the number of pairs $(C',C)$ with $C$ a smooth curve of genus $2$
and $C'$ an \'etale double cover of $C$ with $p$-rank zero.
Our method is to use intersection theory on stacks. In the terminology
of Mumford [M2,~p.293], we intersect the $Q$-classes of $RV_0(\wtM2)$ and
$P^{-1}V_0(\a1)$ inside $R\wtM2$. The intersection has dimension zero
and the $Q$-class of every point is counted with the intersection
multiplicity in the universal deformation space.
This gives the number of points of $V_0(R\wtM2)$ and we need to subtract
the contribution of the reducible curves to obtain
the number of points of $V_0(R\m2)$.

\proclaim (6.1) Theorem. 
The weighted number of isomorphism classes of pairs $(C', C)$ 
with $C$ a smooth curve of genus $2$
and $C'$ an \'etale double cover of $C$ with $p$-rank zero
is given by
$$
\sum_{(C', C)} {{m(C', C)} \over
{|\Aut(C', C)|}} = {1 \over 384} (p-3)(p-1)^2(p+1),
$$
where $m(C',C)\in\Z_{\ge1}$
is the intersection multiplicity in the universal deformation
space and $|\Aut(C', C)|$ denotes the number of
automorphisms of $C$ that fix the point of order $2$ in $J(C)$ 
corresponding to $C'$.

\noindent{\sl Proof.} 
Since $V_0(\a1)$ has class $(p-1)\la1$, the class of $P^{-1}V_0(\a1)$
equals 
$$(p-1)P^*(\la1)=(p-1)\la1'=(p-1)\phi^*\la1,$$ 
by Corollary 3.2.
By [G,~Thm.9.2], the class of $RV_0(\wtM2)$ equals $(p-1)(p^2-1)\phi^*\la2$
(note that $\wtM2=\a2$ as stacks).
To evaluate the degree of the intersection, we apply $\phi_*$. We obtain
$$\phi_*((p-1)\phi^*\la1\cdot (p-1)(p^2-1)\phi^*\la2)
=15(p-1)^2(p^2-1)\la1\la2={1 \over 384}(p-1)^2(p^2-1),
$$
since $\la1\la2={1\over5760}$ on $\mbar2$.
This is the weighted number of points of $V_0(R\wtM2)$.

The number of reducible curves in this intersection can be computed
easily. We need an ordered pair $(E_1,E_2)$ of supersingular elliptic
curves and a point of order $2$ on $E_1$. We find
$$3\cdot{{p-1}\over{24}}\cdot{{p-1}\over{24}}={{1}\over{192}}(p-1)^2.$$
In the next proposition we show that these pairs count with multiplicity 
$p+1$. Therefore the weighted number of points of $V_0(R\m2)$ equals
$$\eqalign{
{1 \over 384}(p-1)^2(p^2-1)-{{1}\over{192}}(p+1)(p-1)^2
&={1 \over 384}(p-1)^2(p+1)((p-1)-2) \cr
&={1 \over 384}(p-3)(p-1)^2(p+1), \cr}
$$
as required.\qed

In \S7 we give a number of examples illustrating this formula. Note
that $V_0(R\m2)$ is empty for $p=3$, but nonempty for $p\ge5$.

\proclaim (6.2) Proposition. The locus $RV_0(\wtM2)$ intersects
 $P^{-1}V_0(\a1)$  with multiplicity $p+1$ 
 at points corresponding 
to stable curves  with two elliptic components. 

\noindent {\sl Proof.} 
For a smooth \'etale double cover $D' \to D$ corresponding to a point
$\eta$ of order $2$ in $J(D)$ the codifferential of
the Prym map $P:R\Mg \to \a{g-1}$ is given by the multiplication map
${\rm Sym}^2(H^0(D,\omega_D \otimes \eta)) \to H^0(D, \omega_D^2)$.
We extend this result to certain \'etale double covers of reducible
curves of genus $2$ in the following manner.

Let $C$ be a stable curve of genus $2$ consisting of two
elliptic components $E_1$ and $E_2$. 
We let ${\rm Def}(C)$ be the space of infinitesimal deformations of $C$. 
Since $\wtM2=\a2$, it can be identified with the space of infinitesimal
deformations of the Jacobian $E_1\times E_2$ of $C$ (as a 
principally polarized abelian surface). We thus get
$$
{\rm Def}(C) \cong {\rm Sym}^2(T_{E_1}\oplus T_{E_2})
= {\rm Sym}^2(T_{E_1}) \oplus (T_{E_1}\otimes T_{E_2})
\oplus {\rm Sym}^2(T_{E_2}),
$$
where $T_{E_i}$ denotes the tangent space to $E_i$ at the origin.
Note that the tangent space and the cotangent space to $\a2$ are
identified via the principal polarization.

Consider now the \'etale double cover $C'$ of $C$ defined by a point 
$\eta$ of order $2$ of $E_1$. 
It consists of the \'etale double cover $E_1'$ of $E_1$ defined by $\eta$
and two copies of $E_2$. Then 
$$H^0(C,\omega_C\otimes \eta)=H^0(E_1,\omega_{E_1}\otimes\eta)\oplus 
H^0(E_2,\omega_{E_2})=H^0(E_2,\omega_{E_2})=T_{E_2}\,.$$
With these identifications, the codifferential of the Prym map
becomes the inclusion
$${\rm Sym}^2(T_{E_2})\to 
{\rm Sym}^2(T_{E_1}) \oplus (T_{E_1}\otimes T_{E_2})\oplus{\rm Sym}^2(T_{E_2})
$$
and the tangent space to the fiber of the Prym map is given by
$${\rm Sym}^2(T_{E_1}) \oplus (T_{E_1}\otimes T_{E_2}).$$

On the other hand, for a $C$ which is supersingular (i.e., both $E_i$ 
are supersingular) one can calculate the directions
in the deformation space ${\rm Def}(C)$
along which supersingularity (or equivalently,
the condition that the $p$-rank is $0$) is preserved. According
to a calculation by Oort [O2,~p.258] we find the $p+1$ directions
$$
\left(\matrix{-\zeta^p & 1  \cr 1  & \zeta \cr }\right) 
\in {\rm Sym}^2(T_{E_1}\oplus T_{E_2}),
$$
where $\zeta^{p+1}=-1$. Here the $i$th coordinate corresponds to $T_{E_i}$.
These tangent vectors are not contained in 
${\rm Sym}^2(T_{E_1}) \oplus (T_{E_1}\otimes T_{E_2})$. \qed

\bigskip
\centerline{\bf \S7. Examples}
\medskip
\noindent
We now give a number of examples to illustrate the formula of
the preceding section. In \S5 we discussed how the $p$-rank of a
hyperelliptic curve may be computed and how the \'etale double
covers of hyperelliptic curves can be described. In the case
where $g(C)=2$ we may assume that $C_2$ is a rational curve
and that $C_1=E$ is an elliptic curve. The Prym variety of $(C',C)$
is isomorphic to $E$. Below, we will identify the \'etale double
cover $C'$ with a point $\eta$ of order $2$ in $J(C)$. With $C$ defined
by $y^2=f(x)$, we give $\eta$ by choosing $2$ of the $6$ branch points.

\medskip\noindent {\bf (7.1) $p=3$.}
Over ${\overline\F}_3$ the elliptic curve with $j=0=1728$
can be written $y^2=x^3-x$. Hence $C$ can be written
$y^2=(x^2+ax+b)(x^3-x)$ with $\eta$ corresponding to $x^2+ax+b$.
Then the Hasse-Witt matrix equals
$$H=\pmatrix{a_2&a_5\cr a_1&a_4\cr}=\pmatrix{-a&1\cr -b&a\cr}$$
and we have the condition $H\cdot H^{(p)}=0$.
In particular, $0=\det(H)=b-a^2$; note that this condition already makes
$C$ singular: $x^2+ax+a^2=(x-a)^2$. So, not only there isn't a smooth $C$
of genus 2 and 3-rank 0 with unramified double cover $C'$ of 3-rank 0
(hence Prym of 3-rank 0), there isn't even a smooth $C$ of 3-rank $\le1$
with Prym of 3-rank 0. 

\medskip\noindent {\bf (7.2) $p=5$.}
Over ${\overline\F}_5$ the elliptic curve with $p$-rank 0
has $j=0$ and can be written $y^2=x^3-1$. Hence $C$ can be written
$y^2=(x^2+ax+b)(x^3-1)$ with $\eta$ corresponding to $x^2+ax+b$.
Now $a(x)=f(x)^2$ and
$$H=\pmatrix{a_4&a_9\cr a_3&a_8}.
$$
One checks that $\det(H)=0$ gives $b=-a^2$ or $b=-2/a$. Since $b=-a^2$
leads to a singular curve ($x^2+ax-a^2=(x-2a)^2$) one is left with $b=-2/a$.
Working out $H\cdot H^{(p)}=0$ leads to the necessary condition
$a^{15}-a^9+1=0$; but $a^3+3=0$ is to be avoided, as it leads to a singular
curve. The residual polynomial in $\F_5[a]$ of degree 12 has 3 irreducible
factors of degree 4. Now note that the automorphism group of the 4-tuple
$\{1,\zeta,\zeta^2,\infty\}$ (with $\zeta$ a primitive cube root of 1)
is acting on the solution set, reflecting the choice of $x$-coordinate in 
writing the elliptic curve with $j=0$ as $y^2=x^3-1$. It is the group $A_4$,
generated by the transformations $x\mapsto\zeta x$ and
$x\mapsto(x+2)/(x-1)$. A calculation shows that the solution set is
one $A_4$-orbit. So there is one isomorphism class of pairs $(C',C)$
with $C'$ of 5-rank 0; it is defined over $\F_5$ and has only
the hyperelliptic automorphism $y\mapsto -y$ (this follows since the
$A_4$-orbit has length 12). 
So the answer to the counting problem is
$${1\over2}={{2\cdot4\cdot4\cdot6}\over{384}} $$
as in the theorem; this proves that the multiplicity $m(C',C)=1$.
Finally, one would of course like to have a model
over $\F_5$. One can try to take $\eta=\{0,\infty\}$ although this is
not guaranteed to work:
$$y^2=x(x^4+b_3x^3+b_2x^2+b_1x+b_0)$$
with the $b_i\in\F_5$. The Prym condition is just $b_2^2+2b_1b_3+2b_0=0$
wheras the condition for $C$ simplifies to $H^2=0$. 
A solution is (for instance):
$$y^2=x(x^4+x^3+2x+3)\qquad\hbox{with}\qquad H=\pmatrix{4&2\cr2&1\cr}. $$

\medskip \noindent {\bf Remark.}
The locus $RV_0(\m2)$ is nonsingular at a point $(C',C)$ with
$a$-number $a(C)=1$, while it has multiplicity $p+1$ at a point
with $a(C)=2$. If the intersection with $P^{-1}V_0(\a1)$ is transversal,
the multiplicity $m(C',C)$ equals $1$ resp.~$p+1$.

\medskip\noindent {\bf (7.3) $p=7$.}
Over ${\overline\F}_7$ we start with the form
$y^2=(x^2+ax+b)(x^3-x)$. To find solutions one can proceed as follows.
The 4 entries of $H\cdot H^{(p)}$ and $\det(H)$ are 5 polynomials
in $a$ and $b$ that have to vanish simultaneously. 
Using Maple, one computes
the resultants modulo 7 w.r.t.~the variable $b$ (resp.~$a$)
of the determinant and each of the 4 entries without problem; taking
the polynomial g.c.d.~modulo 7 of the 4 resultants leads to a 
polynomial in $a$ (resp.~$b$) of degree 65 in both cases. The multiple factors
in the $a$-polynomial are $a^{19}(a+2)^3(a+5)^3$ reducing the number of
possible $a$-values to 43. The factors of the $b$-polynomial all have
multiplicity at least 2, with $(b+6)^{22}(b^2+3b+1)^4b^3$
the factors with multiplicity $>2$, reducing the number of possible $b$-values
to 20. Examining the singular solutions leads to the exclusion of
$a=\pm2$ and of $b=0$, which leaves us with 41 possible $a$-values and
19 possible $b$-values.
Clearly, when $(a,b)$ is a solution, so is $(-a,b)$; more generally,
the automorphism group of the 4-tuple $\{0,1,-1,\infty\}$, the group
$D_4$, is acting on the solution set; the $x$-transformations
$x\mapsto -x$ resp.~$x\mapsto(x-1)/(x+1)$ have the effect
$(a,b)\mapsto(-a,b)$ resp.
$$(a,b)\mapsto\left({{2b-2}\over{b+1+a}},{{b+1-a}\over{b+1+a}}\right).
$$
We examine the orbits of this action. 
It is clear that $(a,b)=(0,1)$ is the unique fixed point; it is in fact
the unique solution with a non-trivial stabilizer.
Here $C$ is the curve $y^2=x^5-x$ which has 48 automorphisms; 
$\eta=\{i,-i\}$; the pair $(C,\eta)$ has 16 automorphisms. Note that
$a(C)=2$ (since $H=0$); hence we expect to count $(C,\eta)$ with
multiplicity $p+1=8$ for a total contribution of $8/16=1/2$.
The pair $(C,\eta)$ is isomorphic to $(C,\{0,\infty\})$.

It is easy to check that the only solution with $a=0$ is $b=1$,
and vice versa. 
For the solutions with $a$ nonzero, 40 possible $a$-values remain. These
fall into 5 orbits. We find only 18 possible
$b$-values since a `ramification' occurs: to $b^2+3b+1=0$
correspond 8 $a$-values (the roots of $(a^4-a^2-1)(a^4-3a^2-1)$)
instead of the expected 4.

So we find 5 more solutions $(C,\eta)$, each having only the automorphisms
$y\mapsto\pm y$; also, $a(C)=1$ for all five. Hence each contributes $1/2$,
so that the answer to the counting problem is
$$6\cdot{1\over2}=3={{4\cdot6\cdot6\cdot8}\over{384}}.
$$
Two of the 5 $D_4$-orbits are given by a polynomial in $\F_7[a]$ of degree 16
(with 4 irreducible factors of degree 4). Hence the 2 corresponding
isomorphism classes $(C,\eta)$ are defined over $\F_{49}$ (a Galois orbit
of length 2). With some work, one finds explicit representatives:
$$
\left\{\eqalign{C&:\quad y^2=x(x^4+x^3+(i+1)x^2+(3i-2)x-i)=
x(x^2+3i-2)(x^2+x-2i+3)\cr
\eta&:\quad x^2+x-2i+3\cr
E&:\quad y^2=x(x^2+3i-2)\qquad(j=1728=-1)\cr}
\right.
$$
with $i=\sqrt{-1}\in\F_{49}$.

The other three $D_4$-orbits are given by a polynomial in $\F_7[a]$ of degree 
24 (with 4 irreducible factors of degree 6). The 3 corresponding
isomorphism classes $(C,\eta)$ are defined over $\F_{343}$ (a Galois orbit
of length 3). One eventually finds explicit representatives:
$$
\left\{\eqalign{C&:\quad y^2=(x^2+1)(x^4+x^3+(\rho^2+\rho+2)x^2
+(3\rho^2-3\rho+1)x+(5\rho^2-3\rho+3))\cr
\eta&:\quad x^2+1 \cr
E&:\quad y^2=x^4+x^3+(\rho^2+\rho+2)x^2
+(3\rho^2-3\rho+1)x+(5\rho^2-3\rho+3) \qquad(j=1728) \cr}
\right.
$$
with $\rho=\root3\of2\in\F_{343}$.

\medskip\noindent {\bf (7.4) $p=11$.}
There are two cases: (A) the Prym variety is the elliptic
curve with $j=1728$; we take $y^2=(x^3-x)(x^2+ax+b)$ as $(C,\eta)$;
(B) the Prym is the elliptic curve with $j=0$; we work with
$y^2=(x^3-1)(x^2+ax+b)$. 

In case (A) one initially finds a polynomial
in $a$ of degree 340; only its 5 linear factors come with (varying)
multiplicities; the reduced polynomial has degree 197.
The initial polynomial in $b$ also has degree 340; every factor
has at least multiplicity 2, and only its 4 linear factors have higher
multiplicities; the reduced polynomial has degree 100. As it turns out,
the linear factors lead either to singular curves or to the 4 solutions
$(a,b)=(0,2)$, (0,6), (3,1), (8,1). These 4 solutions form
a single $D_4$-orbit; hence the corresponding pair $(C,\eta)$ has
4 automorphisms. One easily checks that $a(C)=2$ so that the
expected contribution is $(p+1)/4=12/4=3$.
The remaining $192$ solutions of the reduced polynomial in $a$ form
$24$ orbits; their expected contribution is $24/2=12$, making the
total expected contribution in case (A) equal to $15$.
The Galois orbits have lengths 1, 4, 4, 4, 6, and 6.

In case (B) one initially finds a polynomial
in $a$ of degree 270; only a linear and a quadratic factor come with 
a multiplicity; the reduced polynomial has degree 243.
Exactly the same happens for the polynomial in $b$.
As it turns out, the linear and the quadratic factor
lead to singular curves. 
The remaining $240$ solutions of the reduced polynomial in $a$ form
$20$ $A_4$-orbits; the expected contribution in case (B) is $20/2=10$.
The Galois orbits have lengths 3, 3, 4, and 10.

Hence the total expected contribution equals $25$, as in Theorem 6.1.
So all intersections are transversal.

\medskip\noindent {\bf (7.5) $p=13$.}
Here we find the first example of a non-transversal intersection.
We work with $y^2=(x^3+x+4)(x^2+ax+b)$ with Prym variety the supersingular
elliptic curve $y^2=x^3+x+4$. One initially finds a polynomial in $a$
of degree $466$; a linear and a quadratic factor come with multiplicity $15$
and yield singular curves; the factor $a^4+8a^3+7a^2+8a+9$ comes with
multiplicity $2$; one linear factor doesn't give a solution;
the remaining reduced polynomial has degree $416$. The remaining polynomial
in $b$ has degree $416$ as well. The group $D_2$ acts; if all intersections
were transversal and all points of intersection had $a$-number $1$, then the
total contribution would be $104/2=52$. However, the answer is $105/2$
according to Theorem 6.1. The only possible explanation is that exactly
one of the intersections is non-transversal and has multiplicity $2$.
Indeed, the factor $a^4+8a^3+7a^2+8a+9$ yields that intersection point.
The pair $(C,\eta)$ is defined over $\F_{13}$; an equation is
$$y^2=x^5+x^4+6x^3-2x^2+2x$$
with $\eta=\{0,\infty\}$. This is a nonsingular point of $RV_0(\m2)$
and the tangent line may be given as
$$y^2=x^5+x^4+(t+6)x^3+(2t-2)x^2+(2-t)x$$
with the same $\eta$. The moving elliptic curve is
$$y^2=x^4+x^3+(t+6)x^2+(2t-2)x+2-t$$
and the coefficient of $x^{12}$ in $(x^4+x^3+(t+6)x^2+(2t-2)x+2-t)^6$
is $t^2$ times a unit, which proves the multiplicity statement.

\medskip\noindent {\bf (7.6) $p=17$.}
There are two supersingular $j$-invariants. For the elliptic curve
$y^2=x^3+x+1$ we find a Galois orbit of length $2$ with $a$-number $2$
and $4$ automorphisms, contributing $2(p+1)/4=9$, and a pair $(C,\eta)$
defined over $\F_{17}$ giving a simple tangency and contributing $1$. The
remaining solutions contribute $928/8=116$ for a total contribution of
$126$.

For the elliptic curve $y^2=x^3-1$ we find the first example of a
non-transversal intersection with $a$-number $2$. This proves that
the intersection multiplicity cannot be read off from the
Hasse-Witt matrix. Let the pair $(C,\eta)$ be given by
$$y^2=(x^3-1)(x^2-2x-2).$$
The $18$ tangent directions to $RV_0(\m2)$ can be written down explicitly.
Two of them lie in the tangent plane to the Prym fiber and it turns out
that they yield simple tangencies. The intersection multiplicity equals
$20$. The $A_4$-orbit of this pair has length $3$ and it contributes
$3\cdot20/24=5/2$. There is also an $A_4$-orbit of length $12$
with $a$-number $2$, contributing $12\cdot18/24=9$. The remaining
solutions contribute $732/24=61/2$. The total contribution for
$y^2=x^3-1$ equals $42$.

Adding up the contributions for the two elliptic curves, we obtain
$168$, as in Theorem 6.1.

\bigskip \noindent
{\bf Acknowledgements.} 
We thank Frans Oort for several useful discussions; Proposition 4.3
is due to him. We thank the G\"oran Gustafsson Stiftelse
for support of the second author's visit to KTH in March 2001.
The first author thanks Princeton University and
the Swedish Research Council for support.
Finally, we thank the referee for several useful comments.

\bigskip
\centerline{\bf References}
\medskip
\parindent=0pt

[D1] S. Diaz: A bound on the dimensions of complete subvarieties of 
${\cal M}_g$. Duke Math. J. {\bf 51} (1984), 405--408.

[D2] S. Diaz: Complete subvarieties of the moduli space of smooth curves.  
In {\sl Algebraic geometry, Bowdoin 1985\/} (S. Bloch, ed.), 77--81, 
Proc. Sympos. Pure Math., 46, Part 1, Amer. Math. Soc., Providence, RI, 1987.

[F] C. Faber: A conjectural description of the tautological ring of the 
moduli space of curves. 
In {\sl Moduli of curves and abelian varieties\/} (C. Faber and E. Looijenga, 
eds.), 109--129, Aspects Math., E33, Vieweg, Braunschweig, 1999.

[FL] C. Faber and E. Looijenga: Remarks on moduli of curves.
In {\sl Moduli of curves and abelian varieties\/} (C. Faber and E. Looijenga, 
eds.), 23--45, Aspects Math., E33, Vieweg, Braunschweig, 1999.

[FP] C. Faber and R. Pandharipande: 
Logarithmic series and Hodge integrals in the tautological ring.
With an appendix by Don Zagier.
Michigan Math. J. {\bf 48} (2000), 215--252.

[FC] G. Faltings and C.-L. Chai:
Degeneration of abelian varieties.
Ergebnisse der Mathematik und ihrer Grenzgebiete 22.
Springer-Verlag, Berlin, 1990. xii+316 pp.

[G] G. van der Geer: Cycles on the moduli space of abelian varieties. 
In {\sl Moduli of curves and abelian varieties\/} (C. Faber and E. Looijenga, 
eds.), 65--89, Aspects Math., E33, Vieweg, Braunschweig, 1999.

[GO] G. van der Geer and F. Oort:
Moduli of abelian varieties: a short introduction and survey. 
In {\sl Moduli of curves and abelian varieties\/} (C. Faber and E. Looijenga,
eds.), 1--21, Aspects Math., E33, Vieweg, Braunschweig, 1999.

[HM] J. Harris and D. Mumford:
On the Kodaira dimension of the moduli space of curves.
With an appendix by William Fulton.
Invent. Math.  {\bf 67}  (1982), 23--88.

[Ka] N. Katz: Slope filtration of $F$-crystals.  
Journ\'ees de G\'eom\'etrie Alg\'ebrique de Rennes (Rennes, 1978), Vol. I,  
113--163, Ast\'erisque No. 63 (1979).

[KS] S. Keel and L. Sadun: Oort's conjecture for $\a{g}$.
Preprint 2002, math.AG/0204229, to appear in J. Amer. Math. Soc. 

[Ko] N. Koblitz: 
$p$-adic variation of the zeta-function over families of varieties defined 
over finite fields. Compositio Math. {\bf 31} (1975), 119--218.

[L]  E. Looijenga: On the tautological ring of ${\cal  M}_g$. 
Invent. Math. {\bf 121} (1995),  411--419.

[M-B] L. Moret-Bailly:  Pinceaux de vari\'et\'es ab\'eliennes. Ast\'erisque 
No. 129 (1985). 

[M1] D. Mumford: Prym varieties I. In {\sl Contributions to analysis\/} 
(a collection of papers dedicated to Lipman Bers), 325--350. 
Academic Press, New York, 1974. 

[M2] D. Mumford: Towards an enumerative geometry of the moduli space of curves. 
In {\sl Arithmetic and geometry\/} (M. Artin and J. Tate, eds.), 
Vol. II, 271--328, Progr. Math., 36, Birkh\"auser, Boston, 1983. 

[NO] P. Norman and F. Oort: Moduli of abelian varieties. 
Ann. of Math. {\bf 112} (1980), 413--439.

[O1] F. Oort: Subvarieties of moduli spaces. 
Invent. Math. {\bf 24} (1974), 95--119. 

[O2]  F. Oort: Hyperelliptic supersingular curves. 
In {\sl Arithmetic algebraic geometry (Texel, 1989)\/} 
(G. van der Geer, F. Oort, J. Steenbrink, eds.),
247--284, Progr. Math., 89, Birk\-h\"auser, Boston, 1991.

[O3] F. Oort:  Complete subvarieties of moduli spaces. 
In {\sl Abelian varieties (Egloffstein, 1993)\/} 
(W. Barth, K. Hulek, H. Lange, eds.), 
225--235, de Gruyter, Berlin, 1995.

[SV] K.-O. St\"ohr and J.F. Voloch: 
A formula for the Cartier operator on plane algebraic curves.  
J. Reine Angew. Math.  {\bf 377}  (1987), 49--64.

\bigskip
\noindent Institutionen f\"or Matematik \par
\noindent Kungliga Tekniska H\"ogskolan \par
\noindent S-100 44 Stockholm, Sweden \par
\noindent faber@math.kth.se \par

\bigskip
\noindent Korteweg-de Vries Instituut \par
\noindent Universiteit van Amsterdam \par
\noindent Plantage Muidergracht 24 \par
\noindent NL-1018 TV Amsterdam, The Netherlands \par
\noindent geer@science.uva.nl \par

\bye